\RequirePackage{fix-cm}
\documentclass[english]{scrartcl}
\usepackage[T1]{fontenc}
\usepackage[latin9]{inputenc}
\usepackage{geometry}
\usepackage{amsmath}
\usepackage{amsthm}
\usepackage{amssymb}

\makeatletter
\numberwithin{equation}{section}
\numberwithin{figure}{section}
\@ifundefined{lettrine}{\usepackage{lettrine}}{}

\usepackage{ifxetex}\ifxetex
\usepackage{fontspec}
\else
\fi

\usepackage{fixltx2e}\usepackage{fancyhdr}

\makeatother

\usepackage{babel}
\begin{document}

\title{Epistemological Consequences of the Incompleteness Theorems}

\author{Giuseppe Raguní}

\date{UCAM - Universidad Católica de Murcia, Avenida Jerónimos 135, Guadalupe
30107, Murcia, Spain - graguni@ucam.edu}
\maketitle
\begin{abstract}
After highlighting the cases in which the semantics of a language
cannot be mechanically reproduced (in which case it is called \emph{inherent}),
the main epistemological consequences of the first incompleteness
Theorem for the two fundamental arithmetical theories are shown: the
non-mechanizability for the truths of the first-order arithmetic and
the peculiarities for the model of the second-order arithmetic. Finally,
the common epistemological interpretation of the second incompleteness
Theorem is corrected, proposing the new \emph{Metatheorem of undemonstrability
of internal consistency}.

\end{abstract}

\section{Semantics in the Languages}

Consider an arbitrary language that, as normally, makes use of a countable\footnote{Finite, as the usual alpha-numeric symbols, or, to generalize, infinite-countable.
In this paper we will use either ``countable'' or ``enumerable''
with the same meaning, i. e. to indicate that there exists a biunivocal
correspondence between the considered set and the set of the natural
numbers. } number of characters. Combining these characters in certain ways,
are formed some fundamental strings that we call \emph{terms} of the
language: those collected in a dictionary. When the terms are semantically
interpreted, i. e. a certain meaning is assigned to them, we have
their distinction in adjectives, nouns, verbs, etc. Then, a proper
grammar establishes the rules of formation of sentences. While the
terms are finite, the combinations of grammatically allowed terms
form an infinite-countable amount of possible sentences.

In a non-trivial language, the meaning associated to each term, and
thus to each expression that contains it, is not always unique. The
same sentence can enunciate different things, so representing different
\emph{propositions}. For example, the same sentence ``it is a plain
sailing'' has a different meaning depending on the circumstances:
at board of a ship or in the various cases with figurative sense.
How many meanings can be associated to the same term? That is: how
many different propositions, in general, can we get by a single sentence?
The answer, for any common semantic language, may be amazing.

Suppose we assign to each term a finite number of well-defined meanings.
We could then instruct a computer to consider all the possibilities
of interpretation of each term. The computer, to simplify, may assign
all the different meanings to an equal number of distinct new terms
that it has previously defined. For example, it might define the term
``f-sailing'' for the figurative use of ``sailing'' (supposed
unique). The machine would then be able, using the grammar rules,
to generate all the infinite-countable propositions. In this case
we will say that, in the specific language, the meaning has been \emph{eliminated}\footnote{Just a concise choice rather than \emph{mechanically reproduced}.}\emph{.}
More generally we have this case when the different meanings allowed
for each term are \emph{effectively }enumerable\footnote{A set is called \emph{effectively }enumerable if there is a machine
capable of producing in output all and only its elements.}: even in the case of an infinite-countable amount of meanings, the
computer can define an infinite-countable number of new terms and
associate only one meaning to each term in order to establish a biunivocal
correspondence between sentences and propositions. So, the machine
could list all them by combination.

Hence, by definition, we will say that a language is \emph{inherently
semantic} (i. e. with a non-eliminable meaning) if it uses at least
one term with an amount \emph{not} effectively enumerable of meanings;
with the possibility, which we will comment soon, that this quantity
is even non-enumerable (or \emph{uncountable}). From the fact that
a sentence represents more than one proposition if and only if it
contains at least one term differently interpreted, it follows an
equivalent condition for the inherent semanticity: a language is inherently
semantic if and only if the set of all possible propositions is not
effectively enumerable.

Now, the case of uncountable propositions is really what happens in
every usual natural language. At first, this feature might surprise
or be considered unacceptable: all the meanings that \emph{ever will
be assigned} to any settled word are only a countable number. Indeed,
a \emph{finite} number! But these meanings cannot be specified once
and for all. The fact remains that the \emph{possible} interpretations
of the term vary within an infinite collection; moreover, a collection
not limited by any prefixed cardinality.

Some classic paradoxes can be interpreted as a confirmation of this
property {[}1{]}. The \emph{Richard's} one\footnote{Where firstly it is admitted that all possible semantic definitions
of the real numbers stay in a countable array and, then, one can define,
by a \emph{diagonal argument} {[}2{]}, a real number that is not present
in the array.}, for example, can be interpreted as a meta-proof that the semantic
definitions are not countable, i. e. that they are \emph{conceivably}
able to define each element of a set with cardinality greater than
the enumerable one (and therefore each real number). The proper technique
used in a \emph{diagonal} \emph{Cantorian argument} {[}2{]}, reveals
that the natural language is able to adopt different semantic levels
(or contexts), looking ``after'' (or from the ``outside'') what
was ``before'' (or ``inside'') defined; namely, what was previously
said by the same language. Identical words used in different contexts
have a different meaning and for the number of contexts, including
nested, there is no limit.

On the other hand, the \emph{Berry's }paradox {[}1{]} clearly shows
that a finite amount of symbols, differently interpreted, is able
to define an infinite amount of objects. Here the key of the argument
is again the use of two different contexts to interpret the verb ``define''.

\section{Inherent Semantics in Arithmetics }

Now, on a properly mathematical context, consider the formal first-order
arithmetic theory (\emph{PA}, from \emph{Peano's Arithmetic}). It
is possible to prove that its propositions and theorems are effectively
enumerable {[}3{]}, so \emph{PA} is not inherently semantic. However,
from the first incompleteness Theorem follows that \emph{the truths}
of \emph{PA} are not effectively enumerable. In fact, starting by
\emph{PA}, supposed consistent, add as new axioms the class of all
its \emph{true} statements. We mean: true in the intuitive (or \emph{standard})
\emph{model }(i. e. consistent interpretation), formed by the spontaneous
natural numbers. We have so created a new axiomatic system (\emph{PAT})
that, by construction, is \emph{syntactically complete}\footnote{An axiomatic theory is called \emph{syntactically complete }if, for
any its sentence, it or its negation is a theorem.}: given any sentence, if it is true is a theorem and if it is false
then its negation is a theorem. Now, according to the first incompleteness
Theorem we conclude that the axioms of \emph{PAT}, although enumerable,
cannot be effectively enumerable {[}4-6{]}. Therefore, for what we
said, the truths of the formal first-order arithmetic have to contain
inherent semantics. 

Equivalently, we can say that in the expression ``true statement
in the standard model'' the term \emph{true} has got an amount not
effectively enumerable (although enumerable) of distinct meanings.
So, the phrase belongs to an inherently semantic language.

The interesting epistemological consequence of this fact is that no
machine can be programmed to get all and only the truths of the Arithmetic.
In other words, these truths can only be obtained by the use of not
predetermined criteria, so clarifying, interpreting (or inventing!)
more and more unpredictable meanings for the concept of truth, in
an infinite process. We also can say that the concept of natural number,
though spontaneous and primitive, is not \emph{mechanizable}: its
nature is inherently semantic.

The \emph{PA} theory admits an infinite amount of other interpretations
besides the standard model. The claim to build an arithmetical theory
with the standard interpretation as unique model (briefly: \emph{categorical}),
leads to a more general axiomatic system: the (\emph{full}) second-order
arithmetic (\emph{FSOA}). Here, the induction principle is extended
to \emph{all} the properties of the natural numbers, by a second-order\footnote{For definition, in a \emph{second-order} syntactical expression, the
existential quantifiers $\forall$ and $\exists$ can range not only
over the variable-elements, but also over variable-properties of the
elements of the model. About the dangerous confusions regarding the
expressive order, see {[}7{]}.} axiom\footnote{See for example {[}8{]}, where\emph{ FSOA} is called \emph{AR}.}
Since these properties, as was proved by Cantor, are uncountable {[}2{]},
it is admitted that this theory can express a non-enumerable amount
of propositions. So, in particular, these are not effectively enumerable:\emph{
FSOA}, unlike \emph{PA,} is an inherently semantic discipline.

The categorical nature has, therefore, its price. The (unique) model
of \emph{FSOA} reveals all the peculiarities of the standard, intuitive,
sight of the natural numbers. Specifically, this unique interpretation
has to be able to assign more than a single meaning (indeed a quantity
at least 2\textsuperscript{$\aleph_{0}$}) to at least one sentence,
before to verify the premises of the theory. So, this model cannot
be confused with a conventional uncountable model of a formal, and
maintained formal, axiomatic theory. This latter one, considered for
example for \emph{PA}, will continue to assign \emph{only one} meaning
to each sentence of the system: an uncountable amount of exceeding
elements of the universe will have no representation in \emph{PA}.
Instead, in the\emph{ FSOA }theory the unique model is, in a sense,
``dynamical'' and a proposition, in general, is not fixed only by
its symbolic chain. Semantics is required even to form the propositions!

\section{About an improper interpretation of the second incompleteness Theorem}

The semantical consequences of the second incompleteness Theorem are
often flawed. In reference to a theory that satisfies the same hypothesis
of the first incompleteness Theorem, the second one generalizes the
undecidability to a class of statements which, interpreted in the
\textit{standard model} mean ``this system is consistent''. Its
demonstration, only outlined by Gödel, was published by Hilbert and
Bernays in 1939.

The usual interpretation of this Theorem, object of our criticism,
is that ``every theory that satisfies the hypotheses of the first
incompleteness Theorem cannot prove its own consistency''. It seems
clear, in fact, that the conclusion that a theory cannot prove its
own consistency is valid for \textit{all the classical systems}, including
non-formal! Moreover, this conclusion does not correspond to the second
incompleteness Theorem, but to a new and autonomous metatheorem.

Consider an arbitrary classical axiomatic system. If it is inconsistent,
it is deprived of models and therefore of any reasonable interpretation
of any statement\footnote{More in depth: of any interpretation respecting the principles of
contradiction and excluded third.}. Therefore, only the admitting that a given statement of the theory
means something, implies agreeing consistency. And indubitably this
also applies if the interpretation of the statement is ``this system
is consistent''.

So, if there is no assurance about the consistency of the theory (which,
to want to dig deep enough, applies to any mathematical discipline)
we cannot be certain on any interpretation of its language. For example,
in the case of the usual Geometry, when we prove the pythagorean Theorem,
what we really conclude is ``if the system supports the Euclidean
model (and therefore is consistent), then in every rectangle triangle
c\textsubscript{1}\textsuperscript{2}+c\textsubscript{2}\textsuperscript{2}=I\textsuperscript{2}''.
Certainly, a deduction with an undeniable epistemological worth, still
in the catastrophic possibility of inconsistence.

But now let's see what happens if a certain theorem of a certain theory
is interpreted with the meaning: ``this system is consistent'' in
a given interpretation \textit{M.} Similarly, what we can conclude
by this theorem is really: ``if this system supports the model \textit{M}
(and therefore is consistent), then it is consistent''\textit{. Something
that we already knew and, above all, that does not demonstrate at
all the consistency of the system}\footnote{Just remember that in case of inconsistency we have that every proposition
is a theorem!}\textit{.} Unlike any other statement with a different meaning in
\textit{M}, for this kind of statement we have a peculiar situation:\textit{
bothering to prove it within the theory is epistemologically irrelevant
in the ambit of the interpretation M}. In more simple words, the statement
in question can be a theorem or be undecidable with no difference
for the epistemological view. Just it cannot be the denial of a theorem,
if \textit{M} is really a model. So, in any case, \textit{the problem
of deducing the consistency of the theory is beyond the reach of the
theory itself.} We propose to call \textit{Metatheorem of undemonstrability
internal of consistency} this totally general metamathematic conclusion.

Then, the fact that in a particular and hypothetically consistent
theory, such a statement is a theorem or is undecidable, is depending
on the system and on the specific form of the statement. For theories
that satisfy the assumptions of the first one, the second incompleteness
Theorem guarantees that ``normally'' these statements are undecidable.
``Normally'' just means ``in the standard interpretation'': in
fact the Theorem does not forbid that statements expressing consistency
of the system in different models could be theorems\footnote{Indeed, it seems that these statements really exist {[}9-10{]}.}.
As Lolli says, ``it seems that not even a proof shuts discussions''
{[}11{]}. But in no case this debate can affect the validity of the
proposed Metatheorem.

Summarizing, the second incompleteness Theorem identifies another
class of essentially undecidable statements for any theory that satisfies
the hypothesis of the first one. Whilst the first incompleteness Theorem
determines only the Gödel's statement, the second one extends the
undecidability to a much broader category of propositions. But, contrary
to what is commonly believed, this drastic generalization does not
introduce any new and dramatic epistemological concept about the consistency
of the system. It doesn't so, even if the Theorem were valid \textit{for
every} statement interpretable as ``this system is consistent''
(which, we reaffirm, seems to be false). Because by it in no case
we can conclude that ``the system cannot prove its own consistency'':
this judgment belongs to a different and completely general Metatheorem
which seems never have been stated, despite its obviousness and undeniability\footnote{ The consistency of a theory may be demonstrated only outside the
same, by another external system. For which, in turn, the problem
of consistency arises again. To get out of this endless chain, the
``last'' conclusion of consistency has to be purely metamathematic.
Actually, the most general theory (that demonstrates the consistency
of the ordinary mathematic disciplines) is the formal Set Theory and
the conclusion of its consistency only consists in a ``reasonable
conviction''.}.

Finally we emphasize that the meta-demonstration of the proposed Metatheorem,
since refers to any arbitrary classical system, must consist in a
purely meta-mathematical reasoning (like that one we have presented):
it cannot be formalized.

\end{document}